\documentclass[11pt]{amsart}
\usepackage{amsmath, amsthm, latexsym, amssymb, amsfonts}

\newenvironment{pf}{\proof[\proofname]}{\endproof}
\theoremstyle{plain}
\newtheorem{Th}{Theorem}[section]
\newtheorem{Cor}[Th]{Corollary}
\newtheorem{Prop}[Th]{Proposition}

\newtheorem*{ThA}{Theorem A}
\newtheorem*{ThB}{Theorem B}
\theoremstyle{definition}
\newtheorem{Emp}{}
\newtheorem{Def}[Th]{Definition}
\newtheorem*{Ack}{Acknowledgements}
\theoremstyle{remark}
\newtheorem{Example}{Example}
\newtheorem{Rem}[Th]{Remark}

\newcommand{\cl}[1]{\mathcal{#1}}
\newcommand{\C}{\mathbb C}
\newcommand{\Z}{\mathbb Z}
\newcommand{\R}{\mathbb R}
\newcommand{\G}{\Gamma}
\newcommand{\g}{\gamma}
\newcommand{\D}{\Delta}

\newcommand{\cF}{\cl F}

\newcommand{\cO}{\cl O}

\newcommand{\cX}{\cl X}
\newcommand{\cY}{\cl Y}
\newcommand{\T}{{\mathbb T}}

\newcommand{\ov}{\overline}

\newcommand{\p}{\partial}
\newcommand{\sig}{\sigma}
\newcommand{\Sig}{\Sigma}

\newcommand{\res}{\operatorname{res}}

\newcommand{\tres}{\operatorname{res}^\T}

\newcommand{\sgn}{\operatorname{sgn}}

\newcommand{\St}{\operatorname{St}}
\newcommand{\Spec}{\operatorname{Spec}}
\newcommand{\Hom}{\operatorname{Hom}}

\newcommand{\symb}{\langle f_0,\dots,f_{n}\rangle}
\newcommand{\tsymb}{[f_0,\dots,f_{n}]}

\newcommand{\rs}[1]{Section~\ref{S:#1}}

\newcommand{\rp}[1]{Proposition~\ref{P:#1}}
\newcommand{\rr}[1]{Remark~\ref{R:#1}}
\newcommand{\rex}[1]{Example~\ref{ex:#1}}
\newcommand{\re}[1]{(\ref{e:#1})}
\newcommand{\rco}[1]{Corollary~\ref{C:#1}}
\newcommand{\rt}[1] {Theorem~\ref{T:#1}}
\newcommand{\rd}[1]{Definition~\ref{D:#1}}

\begin{document}


\title{Residues and tame symbols on toroidal varieties}
\author{Ivan Soprounov}
\email{isoprou@math.umass.edu}
\address{Department of Mathematics and Statistics\\
University of Massachusetts, Amherst\\ Amherst, MA 01003\\ USA}

\begin{abstract} We introduce a new approach to
the study of a system of algebraic equations in $({\C}^\times)^n$ whose
Newton polytopes have sufficiently general relative positions.
Our method is based on the theory of Parshin's
residues and tame symbols on toroidal varieties. It provides a
uniform algebraic explanation of the
recent result of Khovanskii on the product
of the roots of such systems and the Gel'fond--Khovanskii result on
the sum of the values of a Laurent polynomial over the roots of such
systems, and extends them to the case of an algebraically
closed field of arbitrary characteristic.
\end{abstract}

\maketitle


\contentsline {section}{\tocsection {}{1}{Introduction}}{1}
\contentsline {section}{\tocsection {}{2}{Degree of Polyhedral Maps and Combinatorial Coefficient}}{4}
\contentsline {section}{\tocsection {}{3}{Toroidal Symbol}}{7}
\contentsline {section}{\tocsection {}{4}{Toroidal Residue}}{15}
\contentsline {section}{\tocsection {}{5}{Application to Systems of Equations in the Torus}}{19}
\contentsline {section}{\tocsection {}{6}{Appendix: Parshin's Reciprocity Laws}}{22}
\contentsline {section}{\tocsection {}{}{References}}{25}

\section{Introduction}

\begin{Emp}
The classical residue formula says that the sum of the residues
of a rational 1-form $\omega$ over all points of a complex projective curve $X$
is zero:
$$\sum_{x\in X}\res_x\omega=0.$$
The standard proof of this formula uses Stokes theorem. In
``Algebraic groups and class fields'' J.-P.~Serre gives a purely
algebraic proof of the residue formula which works over any algebraically
closed field even of positive characteristic \cite{S}.

In class field theory
the residue formula has a multiplicative cousin --- Weil's
reciprocity: the product of the tame symbols of any two rational
functions $f_0$, $f_1$ over all points of a projective curve $X$ is one
\cite{S}:
$$\prod_{x\in X}\langle f_0,f_1\rangle_x =1.$$

In the 1970's A.~Parshin constructed higher-dimensional class
field theory where he generalized the residue and the tame symbol \cite{P1}.
Given an $n$-dimensional algebraic variety $X$ and a rational form
$\omega$ on $X$, Parshin defines the residue $\res_F\omega$ at each
complete flag $F: X_0\subset X_1\subset\dots\subset X_{n-1}\subset
X$ of irreducible subvarieties of $X$. Similarly, given any $n+1$
rational functions $f_0,\dots,f_{n}$ on $X$ he defines the tame
symbol $\symb_F$ at each such flag $F$. Parshin's residue and
symbol satisfy not one but many reciprocity laws: Fix all
subvarieties in the flag $F$ except one, say $X_i$. Then the sum
of the residues (product of the symbols) over all possible
irreducible subvarieties $X_i$ that can appear in the $i$-th slot
of $F$ is zero (one) (\rt{Parshinrecip} in the appendix).

The aim of the present paper is to look at some recent
results of the theory of Newton polytopes from
the point of view of this general theory of Parshin.
More specifically, consider a system of Laurent polynomial equations
in the $n$-torus:
\begin{equation}\label{e:system0}
f_1(t)=\dots=f_n(t)=0,\quad t\in(\C^\times)^n.
\end{equation}
Suppose the Newton polytopes $\D_1,\dots,\D_n$ of the $f_i$ have
sufficiently general relative positions (see \rd{developed}). 
Then the system has a finite number of roots. O.~Gel'fond and A.~Khovanskii
 proved the following result \cite{G-Kh}.

\begin{ThA} The sum of the values of a Laurent polynomial
$f_0$ over the roots of \re{system0} counting
multiplicities is equal to
$$\sum_{A}(-1)^nc(A)\res_{A}
\left(f_0\frac{df_1}{f_1}\wedge\dots\wedge\frac{df_n}{f_n}\right),$$
where the sum is over the vertices $A$ of $\D=\D_1+\dots+\D_n$,
$\res_{A}\left(f_0\frac{df_1}{f_1}\wedge\dots\wedge\frac{df_n}{f_n}\right)$
is the {\it residue at a vertex} (an explicit rational function in
the coefficients of the $f_i$), and $c(A)$ is the {\it
combinatorial coefficient} (an integer that reflects the
combinatorial structure of the polytope $\D$ near the vertex~$A$).
\end{ThA}

This is, in fact, a particular case of their {residue formula}
\cite{G-Kh2} for the sum of
the Grothendieck residues of a rational $n$-form which is regular
in $(\C^\times)^n\setminus\{f_1\dots f_n=0\}$,
over the roots of the system~\re{system0}. The proof of the
residue formula is topological and uses toric
compactifications.

The following is a generalized Vieta formula for the product of
roots of \re{system0} obtained by Khovanskii \cite{Kh}.

\begin{ThB} The product of the values of a Laurent monomial $f_0$
over the roots of \re{system0} counting multiplicities is equal to
$$\prod_{A}[f_0,\dots,f_n]_A^{(-1)^nc(A)},$$
where the product is over the vertices $A$ of
$\D=\D_1+\dots+\D_n$, $[f_0,\dots,f_n]_A$ is the {\it symbol
at a vertex} (an explicit Laurent monomial in the coefficients of the
$f_i$), and $c(A)$ is the combinatorial coefficient.
\end{ThB}

The proof of this theorem uses the polyhedral homotopy method and
regular subdivisions of polytopes.

The relation between Theorem A and Theorem B appears to be the same as the
one between the 1-dimensional residue formula and Weil's
reciprocity. Moreover, the number $[f_0,\dots,f_n]_A$ is
defined similar to Parshin's tame symbol. This gives a motivation
to search for a uniform explanation of these results in terms of
the theory of residues and tame symbols.

The main obstruction to this is the notion of combinatorial coefficients
$c(A)$ since they are defined as the local degrees of certain real
non-algebraic maps. In the present paper we give an explicit
algebraic description for the combinatorial coefficients as a
signed number of certain complete flags of faces of $\D$, thus
putting them in the framework of Parshin's theory. (Similar
description of the combinatorial coefficient
was obtained by O. Gel'fond \cite{G} for some special collections
of polytopes.) We provide
a uniform algebraic proof of Theorems A and
B based on Parshin's theory for toroidal varieties.
We also extend them to the case of an arbitrary
algebraically closed field.
\end{Emp}

\begin{Emp} The material of the paper is organized as follows. In
section 2 we give an explicit formula for the degree of a map
of polyhedral sets defined by some combinatorial data. As an
application we obtain a new formula for the combinatorial
coefficient.

In sections 3 and 4 we consider Parshin's theory for toroidal
pairs. A toroidal pair $(X,D)$ consists of a normal variety $X$
and a codimension 1 subvariety $D$ such that locally at each point,
$X$ is analytically isomorphic to an affine toric variety $X_\sig$
where the branches of $D$ correspond to the invariant divisors of
$X_\sig$. We define residue and tame symbol at each point $x\in X$
for which the corresponding affine variety $X_\sig$ has a
0-dimensional orbit. This generalizes the notions of
residue and symbol at a vertex from Theorems A and B.
Our definition is similar to the one of Parshin,
but does not involve any particular choice of a complete flag.
Using the algebraic description of the
combinatorial coefficients we prove general results about symbol
and residue on toroidal pairs (\rt{mainsymb}, \rt{mainres}). In
section 5 we show how these results imply Theorems A and B
for arbitrary algebraically closed fields.

Finally, in the appendix we include the definition of Parshin's
residue and tame symbol and formulate the higher-dimensional
reciprocity laws.
\end{Emp}

\begin{Emp}{\it Remarks.} There is a topological
construction for the tame symbol based on Deligne's proof of
Weil's reciprocity (see \cite{BM}). This construction provides
a uniform topological explanation of the product of roots formula
and the residue formula. Our approach is algebraic and
works over algebraically closed fields of arbitrary
characteristic.

A different formula for the product of the roots of a system \re{system0} 
can be derived from Poisson's formula 
for the mixed resultant due to Pedersen and Sturmfels \cite{PS}. 
In this formula the product of the values of a monomial $f_0=ct^m$ over
the roots is represented by the product of the mixed resultant
of $f_0,\dots,f_n$ (which in this case is just $c$ to the power
of the mixed volume of $\D_1,\dots,\D_n$) and the facet resultants 
to certain powers. 
The assumption on the Newton polytopes of \re{system0} implies that
each facet resultant is a monomial in one of the coefficients
of the system. To give an idea of how this formula is related to the one
in Theorem B let us assume that $f_0=c$, for $c\neq 0,1$. 
Then Poisson's formula gives the inductive formula for the mixed 
volume (e.g. see \cite{BZ}, p.166), 
whereas Theorem B gives Khovanskii's formula for the mixed volume 
in terms of combinatorial coefficients \cite{Kh}. 

Parshin's residue is closely related to the
toric residue defined by D. Cox \cite{C}. For different applications
of residues in toric geometry we refer the reader to works of E.~Cattani,
D.~Cox, A.~Dickenstein, and B.~Sturmfels \cite{C-C-D, C-D, C-D-S}.
\end{Emp}

\begin{Emp}
In the paper $k$ is always an algebraically closed field. A
variety is a reduced separated scheme of finite type over $k$, a
subvariety is a reduced subscheme. By $\T$ we denote the algebraic
$n$-torus over $k$, $\T=(k^\times)^n$, and
$M=\Hom_{alg.\/gp}(\T,k^\times)$ the abelian rank $n$ group of
characters of~$\T$. Finally, $X_\sig$ denotes the affine toric variety
$\Spec k[\sig\cap M]$ defined by a convex rational polyhedral
cone $\sig$ in $M\otimes\R$.
\end{Emp}

\begin{Ack} The results of the present paper constitute a
part of the author's Ph.D. thesis \cite{T}. I am grateful
to Askold Khovanskii for stating the problem, his constant support
and numerous discussions.
\end{Ack}

\section{Degree of Polyhedral Maps and Combinatorial Coefficient}\label{S:combcoeff}

In this section we show how to compute the degree of a map
between two polyhedral sets which is defined by a map of the partially
ordered sets of their faces. As an application we obtain an
explicit combinatorial formula for the combinatorial coefficient.

\subsection{Polyhedral maps}
A polyhedral set is a finite union of convex compact polytopes
intersecting in faces. We will assume that all the polytopes are
embedded in a Euclidean space $E$ of some big dimension. Then a
polyhedral set is a topological space with the topology inherited
from $E$. The dimension of a polyhedral set is the maximum of dimensions
of the polytopes it contains. A polyhedral set is oriented if
every polytope it contains is oriented.

Let $X$ be a polyhedral set and $\cF(X)$ the set of all faces of all polytopes appearing
in $X$. The set $\cF(X)$ is a finite partially ordered set by
inclusion.

Consider two polyhedral sets $X$ and $Y$, and fix a map
$\psi:\cF(X)\to\cF(Y)$ that preserves the partial ordering. A
continuous piecewise linear map $f_\psi:X\to Y$ is called a
{\it polyhedral map associated with $\psi$}
if $f_\psi(G) \subset \psi(G)$ for every face $G\in\cF(X)$.

Given any map $\psi:\cF(X)\to\cF(Y)$ there exists polyhedral
map $f_\psi:X\to Y$. Moreover all such maps are homotopy
equivalent within the class of all polyhedral maps associated with $\psi$.

\begin{Prop}\label{P:homotopy}
Each map $\psi:\cF(X)\to\cF(Y)$ that respects the partial
ordering defines a homotopy class of polyhedral maps $f_\psi:X\to Y$
associated with $\psi$.
\end{Prop}

\begin{pf} First, for every map of partially ordered sets
$\psi:\cF(X)\to\cF(Y)$ we construct a continuous piecewise
linear map $f_\psi:X\to Y$ as follows.

Fix barycentric subdivisions of $X$ and $Y$. Note that for any
polyhedral set $X$ there is a one-to-one correspondence between the
set of all simplices in a barycentric subdivision of $X$ and the
set of all chains in $\cF(X)$. Consider a $k$-simplex
$\D^k$ in the subdivision of $X$. It corresponds to a chain
$X_0\subset\dots\subset X_k$ in $\cF(X)$. Let
$\psi(X_0)\subseteq\dots\subseteq\psi(X_k)$ be its image. It
corresponds to a unique simplex (possibly of smaller dimension)
in the subdivision of $Y$ which we denote by $\psi(\D^k)$.
Since there is a unique linear map
between two simplices that maps vertices of one simplex to the
prescribed vertices of the other simplex, we get a map $f_\psi:X\to
Y$ that sends each simplex $\D^k$ to the corresponding simplex
$\psi(\D^k)$. Clearly this map agrees on the common faces of
simplices of the subdivision and, hence, is continuous piecewise
linear. By construction, $f_\psi(G)\subset\psi(G)$ for any $G\in\cF(X)$.

Now suppose $f_\psi$ and $f'_\psi$ are two polyhedral maps associated 
with $\psi$.
Then for each $0\leq t\leq 1$ the map $f^t_\psi=(1-t)f_\psi+tf'_\psi$
is also associated with $\psi$. Indeed, fix a face $G\in\cF(X)$. Then
every point $x\in G$ is mapped to a point $f^t_\psi(x)$ on the segment joining
$f_\psi(x)$ and $f'_\psi(x)$. Since both $f_\psi(x)$ and $f'_\psi(x)$ belong to
the face $\psi(G)$, $f^t_\psi(x)$ also does. Therefore 
$f^t_\psi(G)\subset\psi(G)$.
\end{pf}


\subsection{Flags and degree of polyhedral maps}

Consider an oriented polyhedral set $X$. Let $\cX:\, X_0\subset\dots\subset X_{n}$, $\dim X_i=i$ be a complete flag
in $X$, i.e. a maximal chain of elements of $\cF(X)$. With the flag $\cX$
we associate an ordered set of vectors $(e_1,\dots,e_{n})$, where $e_i$ begins at
$X_0$ and points strictly inside~$X_i$. Define the {\it sign} of $\cX$ to
be~$1$ if $(e_1,\dots,e_{n})$ gives a positive oriented
frame for $X_{n}$ and $-1$, otherwise. It is easy to
see that the sign does not depend on the choice of vectors
$e_1,\dots,e_{n}$. We denote it by $\sgn\cX$.

Now consider a map of partially ordered sets
 $\psi:\cF(X)\to\cF(Y)$, where $X$ and $Y$ are $n$-dimensional oriented polyhedral sets.
Every polyhedral map $f_\psi:X\to Y$ associated with $\psi$ induces a map of the $n$-th homology groups:
$$H_n(f_\psi): H_n(X)\to H_n(Y).$$
By \rp{homotopy} this map is the same for all choices of $f_\psi$. We call
it the {\it degree} map of $\psi$.
We will be concerned with the case when both groups $H_n(X)$ and $H_n(Y)$ are
isomorphic to $\Z$. (This is true, for example, when $X$ and $Y$ are the boundaries
of ($n+1$)-dimensional polytopes.) Then the degree map is the multiplication
by an integer, which we denote by $\deg(\psi)$.
In the next theorem we show how to compute $\deg(\psi)$ 
as a signed number of certain complete flags in $X$.

Let $\cX:\, X_0\subset\dots\subset X_n$ and $\cY:\, Y_0\subset\dots\subset Y_n$
be complete flags in $X$ and $Y$, respectively. We will write
$\psi(\cX)=\cY$ if and only if $\psi(X_i)=Y_i$ for all $0\leq i\leq n$.
Define the {\it preimage} of $\cY$ under $\psi$ to be the set of all $\cX$
such that $\psi(\cX)=\cY$.

\begin{Th}\label{T:flags} Let $X$ and $Y$ be two polyhedral sets
as above, and $\psi:\cF(X)\to\cF(Y)$ a map of partially
ordered sets of their faces. Fix any complete flag $\cY$ in $Y$.
Then the degree of $\psi$ is equal to the sign of $\cY$ times the signed
number of all complete flags $\cX$ in $X$ in the
preimage of $\cY$ under $\psi$:
$$\deg(\psi)=\sgn\cY\!\sum_{\psi(\cX)=\cY}\sgn\cX.$$
\end{Th}

\begin{pf}
By \rp{homotopy} we can choose any function in the homotopy class
defined by $\psi$. We take $f_\psi$ to be the piecewise linear function
constructed in the proof of \rp{homotopy} using barycentric subdivisions
of $X$ and $Y$. We view $f_\psi$ as a simplicial map between two simplicial complexes.

Fix any positive oriented $n$-simplex $\D_Y^n$ in the barycentric
subdivision of $Y$. Then the degree of $f_\psi$ is the number of all
$n$-simplices $\D_X^n$ in $X$ that are mapped to $\D_Y^n$; each
simplex being counted with either sign plus or sign minus according to its
orientation. Recall that the $n$-simplex $\D_Y^n$ corresponds to a complete
flag $\cY$ in $Y$ and every $n$-simplex $\D_X^n$ corresponds to some complete
flag $\cX$ in $X$. Clearly, the orientation of $\D_X^n$ coincides with the
sign of the corresponding flag $\cX$, and $f_\psi(\D_X^n)=\D_Y^n$ if
and only if $\psi(\cX)=\cY$. It remains to notice that if we fix a negative
oriented $\D_Y^n$ the number we obtain is the negative degree of $f_\psi$.
\end{pf}


\subsection{Combinatorial coefficient}\label{S:cc}

The combinatorial coefficient is a local analog of the
degree considered above.

Let $\sig\subset\R^n$ be a convex polyhedral $n$-dimensional cone in
with apex~$A$. Consider an ordered collection $D=(D_1,\dots,D_m)$ 
of $m$ distinct non-empty closed subsets 
of $\sig$, where $m\leq n$ and each set $D_i$ is a union of facets of $\sig$.
Assume that they cover the boundary of $\sig$ and if $m=n$ the
apex $A$ is the only face of $\sig$ which is covered by all of
them:
\begin{equation}\label{e:covercone}
\partial\sig=D_1\cup\dots\cup D_m,
\quad\text{ if }\ m=n\ \text{ then }\ D_1\cap\dots\cap
D_{n}=\{A\}.
\end{equation}

A continuous map $g:\sig\to\R^n$ is called a {\it characteristic
map} of the covering \re{covercone} if for each $1\leq i\leq n$
the $i$-th component $g_i$ of $g$ is non-negative and vanishes
precisely on those faces of $\sig$ that belong to $D_i$. It is
easy to see that all characteristic maps send the boundary of
$\sig$ to the boundary of the positive octant $\R^n_+$ such that
$g^{-1}(0)\subseteq\{A\}$, and they are homotopy equivalent within
the class of such maps.

\begin{Def}\label{D:covering}
The local degree of the germ at $A$ of the restriction of a
characteristic map to the boundary of $\sig$, $$\bar
g:(\partial\sig,A)\to(\partial\R^n_+,0)$$ is called the {\it
combinatorial coefficient} of the covering \re{covercone}.
\end{Def}

Clearly, the combinatorial coefficient is zero unless $m=n$. In
the case when $m=n$ \rt{flags} provides us with a description of
the combinatorial coefficient as the number of certain complete
flags of faces of $\sig$, counted with signs.

For a cone $\sig\subset\R^n$ we let $\cF(\p\sig)$ denote the partially
ordered by inclusion set of the proper faces of $\sig$. With a covering
\re{covercone} we associate a map $\phi:\cF(\p\sig)\to\cF(\p\R^n_+)$
by putting
$$\phi(\tau)=\R^n_+\cap\{y_{i_1}=\dots=y_{i_k}=0\}$$
if and only if $\tau$ is a common face of $D_{i_1},\dots D_{i_k}$ for 
$1\leq i_l\leq n$, and $k$ is maximal. Here $(y_1,\dots,y_n)$ 
is a coordinate system for $\R^n$.

For any complete flag $\g_0\subset\dots\subset\g_{n-1}$ in
$\p\R^n_+$ define its {\it preimage} under $\phi$ as the
set of all complete flags $\sig_0\subset\dots\subset\sig_{n-1}$ 
in $\p\sig$ such that $\phi(\sig_i)=\g_i$, $0\leq i\leq n-1$. 
Note that the preimage of any flag under $\phi$ is empty if $m<n$.

Fix the standard orientation of $\R^n$. We orient the boundary of
every $n$-dimensional cone in $\R^n$ in accordance
with this fixed orientation. As before, define the {\it sign} of a 
complete flag $\sig_0\subset\dots\subset\sig_{n-1}$
to be $1$ if it gives a positive oriented frame for $\sig_{n-1}$, and
$-1$ otherwise.

\begin{Th}\label{T:cccone} The combinatorial coefficient of a
covering \re{covercone} is equal to the signed number of all complete flags
in the preimage of any positive complete flag under~$\phi$.

In particular, if $m=n$ the combinatorial coefficient is equal to
the signed number of all complete flags
$\sig_0\subset\dots\subset\sig_{n-1}\subset\sig$,
where $\sig_i$ is a common face of
$D_{i+1},\dots,D_{n}$ of dimension~$i$.
\end{Th}

\begin{pf} The case $m<n$ is obvious, so we assume that $m=n$.
To be able to apply \rt{flags} we ``compactify''
the cones $\sig$ and $\R^n_+$. Consider a pyramid
with the vertex $A$ and base $D_0$ which is a cross
section of $\sig$ by a generic hyperplane.
Let $X$ be the oriented boundary of the pyramid.
Next consider the standard $n$-dimensional simplex defined in
$\R^{n}$ as the convex hull of the origin and the endpoints
of the standard basis vectors. Let $Y$ be its oriented boundary.

The subsets $D_1,\dots,D_n$ of $\sig$ along with the base $D_0$
form a covering of $X$. Define the map $\psi:\cF(X)\to\cF(Y)$,
by putting $\psi(G)=\{y_{i_1}=\dots=y_{i_k}=0\}$ if and only if $G$
is a common face of $D_{i_1},\dots, D_{i_k}$ for $0\leq
i_l\leq n$, and $k$ is maximal. Here $(y_0,y_1\dots,y_n)$ are the barycentric
coordinates for the simplex.

Note that the restriction of a characteristic map $g:\sig\to\R^n$
to $X$ defines a polyhedral map $f_\psi:X\to Y$ associated with $\psi$.
According to \rt{flags}, the degree of $\psi$ is equal to the
number of complete flags of faces of $X$ counted with signs in the
preimage under $\psi$ of any positive complete flag of faces of $Y$.
For example one can take the flag
$$\{y_1=\dots=y_n=0\}\subset\{y_2=\dots=y_{n}=0\}\subset\dots
\subset\{y_n=0\}.$$
\end{pf}

\begin{Rem}\label{R:comborder}
Notice that since there are $n!$ complete flags
in $\R^n_+$ we obtain $n!$ formulae for the combinatorial
coefficient. If $m=n$ a choice of a complete flag corresponds
to an order of $D_1,\dots,D_n$, thus we can say that the
combinatorial coefficient is skew-symmetric in $D_1,\dots,D_n$.
\end{Rem}


\section{Toroidal Symbol}\label{S:symbol}

The toroidal symbol is a slight modification of Parshin's tame symbol for
toroidal varieties. More precisely, consider a pair $(X,D)$ consisting
of a normal variety $X$ and a codimension 1 subset $D$ such that
in a formal neighborhood of each point, $X$ is isomorphic to an
affine toric variety $X_\sig$ and $D$ corresponds to the invariant
divisor $X_\sig\setminus\T$. We distinguish special points on
$X$ for which the corresponding toric variety has a 0-dimensional
orbit. At each such point $x\in X$ the toroidal symbol associates a
non-zero element $\tsymb_x$ of the base field to every collection of $n+1$
rational functions $f_0,\dots,f_{n}$ on $X$ with divisors in $D$.

Suppose the irreducible components of $D$ are divided into $2n$
groups $D_1',\dots,D_n'$, $D_1'',\dots,D_n''$ (where $n=\dim X$) and assume that
the sets $S'=D_1'\cap\dots\cap D_n'$ and $S''=D_1''\cap\dots\cap
D_n''$ consist of special points only. The main result of this
section is a certain reciprocity between the products of symbols
over $S'$ and $S''$.


\subsection{Toroidal pair}\label{S:toroidalpair}

Here we recall the definition of a toroidal pair. A detailed treatment
of toroidal pairs is given in \cite{TrEm} where they are
called toroidal embeddings without self-intersections. We use Danilov's
terminology from \cite{D}.

Let $X$ be a normal
$n$-dimensional variety over an algebraically closed field $k$.
Let $D$ be a closed subset of $X$ every irreducible component of
which is a codimension 1 normal subvariety of $X$. We say that the
pair $(X,D)$ is {\it toroidal} if for every closed point $x\in X$ there
exists an $n$-dimensional algebraic torus $\T$, an affine toric variety
$X_\sig$, corresponding to a rational convex
$n$-dimensional cone
$\sig$, and a point $x_0$ in $X_\sig$, such that $(X,D,x)$ is
formally locally isomorphic to $(X_{\sig},X_\sig\setminus\T,x_0)$.
The latter means that there exists an isomorphism of the formal completions of the
local rings
$$\widehat\cO_{X,x}\cong\widehat\cO_{X_{\sig},x_0},$$
such that the image of the ideal of $D$ is mapped to the image of
the ideal of $X_\sig\setminus\T$. We call $(X_\sig,x_0)$ a {\it
local model} of $(X,D)$ at $x$.

\smallskip

Consider the $n$-form
$$\omega_0=\frac{dt_1}{t_1}\wedge\dots\wedge\frac{dt_n}{t_n},$$ where
$(t_1,\dots,t_n)$ are coordinates in $\T$. Automorphisms of $\T$
correspond to monomial changes of coordinates
\begin{equation}\label{e:monomialchange}
u_i=t_1^{q_{i1}}\dots t_n^{q_{in}},\ \ 1\leq i\leq n,\quad Q=(q_{ij})\in
GL(n,\Z).
\end{equation}
We will write $u=t^Q$ to denote the monomial change of coordinates
\re{monomialchange}.
Note that the form $\omega_0$ is preserved under monomial changes
of coordinates  with $\det Q=1$, and changes the
sign when $\det Q=-1$. Therefore, $\omega_0$ provides an analog of
orientation on $X_\sig$.

Furthermore, the choice of coordinates in $\T$ defines an
orientation of the space of characters $M\otimes\R$. Monomial
changes of coordinates \re{monomialchange} preserve this
orientation if and only if $\det Q=1$. Therefore, the orientation
of $M\otimes\R$ and hence of $\sig$ is uniquely determined by the
form $\omega_0$.

We call $(X_\sig,x_0,\omega_0)$ an {\it equipped local model of
$(X,D)$ at~$x$} assuming that the form $\omega_0$ is fixed and the cone
$\sig$ is oriented accordingly.

\smallskip

Let $D=\bigcup_{i\in I}E_i$ be the decomposition of $D$ into
irreducible components. The components of the sets $\bigcap_{i\in
J} E_i\setminus\bigcup_{i\not\in J}E_i$ (where $J\subset I$) are
non-singular and define a stratification of~$X$ (see \cite{TrEm}, p.
57). In particular, $X\setminus D$ is non-singular. The components
of $\bigcap_{i\in J} E_i$ are normal and are the closures of the
strata. Furthermore, for each $x\in X$ the closures of strata
which contain $x$ correspond formally to the closures of the
orbits in a local model $(X_{\sig},x_0)$ at $x$.

We denote by $\St_i(X)$ the set of all $i$-dimensional strata, and
by $\ov\St_i(X)$ the set of the closures of the $i$-dimensional
strata. Note that if $x\in\St_0(X)$ then in every local model
$(X_\sig,x_0)$ at $x$ the cone $\sig$ has an apex and $x_0$ is the
closed orbit in $X_\sig$.

In the next proposition we describe what coordinate
transformations relate different local models at a point
$x\in\St_0(X)$.

\begin{Prop}\label{P:localmodels}
Let $(X,D)$ be a toroidal pair, $x\in\St_0(X)$. Then for any two
local models $(X_\sig,x_0)$ and $(X_{\sig'},x_0')$ at $x$, every
isomorphism
$$\pi:\widehat\cO_{X_\sig,x_0}\cong\widehat\cO_{X_{\sig'},x_0'}$$
that maps the image of the ideal of $X_\sig\!\setminus\T$ to the
image of the ideal of $X_{\sig'}\!\setminus{\T}'$, is induced by a
change of coordinates of the form
$$u_i=\phi_it_1^{q_{i1}}\dots t_n^{q_{in}},
\quad \phi_i\in\widehat\cO_{X_{\sig},x_0}^{\times},\ \ 1\leq i\leq n,\ \
Q=(q_{ij})\in GL(n,\Z),$$ where $t_1,\dots,t_n$ and
$u_1,\dots,u_n$ are coordinate functions on the tori $\T$ and
${\T}'$, respectively.
\end{Prop}

\begin{pf} Let $\Sig(x)$ be the union of all strata $Z$
whose closure $\overline Z$ contains $x$. Denote by $M(x)$ the
group of the Cartier divisors on $\Sig(x)$, supported on
$\Sig(x)\cap D$, and by $M(x)_+$ the subsemigroup  of effective
divisors. For each local model $(X_\sig,x_0)$ at $x$, $M(x)$ is
canonically isomorphic to the group of characters $M$ of $X_\sig$,
and $M(x)_+$ is canonically isomorphic to the semigroup $\sig\cap
M$ (see \cite{TrEm}, p. 61). Therefore the semigroups $\sig\cap M$ and
$\sig'\cap M'$ are isomorphic. In coordinates
$t_1,\dots,t_n$ and $u_1,\dots,u_n$ it corresponds to a monomial
transformation $u=t^Q$, for $Q\in GL(n,\Z)$.

To describe all isomorphisms
$\pi:\widehat\cO_{X_\sig,x_0}\cong\widehat\cO_{X_{\sig'},x_0'}$ it
suffices to describe all automorphisms $\alpha$ of
$\widehat\cO_{X_\sig,x_0}$ that fix the orbits of $X_\sig$. Let
$t_1,\dots,t_n$ be coordinates in~$\T$. Then the ring
$\widehat\cO_{X_\sig,x_0}$ can be identified with the ring of all
formal power series in $t_1,\dots,t_n$ supported
in $\sig\cap M$, where $M$ is identified with $\Z^n$. Denote this
ring  by $A$. Let~$S$ be a multiplicative subset of $A$ consisting
of all elements $\phi t^a$, where $a\in\sig\cap M$ and $\phi$ is
an invertible element of $A$. Then for every automorphism
$\alpha$, we have $\alpha(S)\subseteq S$. Indeed, since $\alpha$ fixes the
orbits of $X_\sig$, it maps every ideal $(t^a)$ to itself. Thus
$\alpha(t^a)=\phi t^a$, for some invertible $\phi$. Therefore,
$\alpha$ induces an automorphism $\alpha_S$ of the localization
$A_S$. Note that $t_1,\dots,t_n\in A_S$, since the elements of
$\sig\cap M$ generate $M$ as a group. Therefore, for each $1\leq
i\leq n$, $\alpha_S(t_i)=\phi_it_i$ for some invertible $\phi_i$.

Conversely, every map $t_i\mapsto \phi_it_i$, $1\leq i\leq n$,
$\phi_i\in A^\times$ defines an automorphism $\alpha$ of $A$,
which fixes the orbits. Indeed, for every element $f\in A$,
$f=\sum_{a}\lambda_a t^a$, $a\in\sig\cap M$, put
\begin{equation}\label{e:automorphism}
\alpha(f)=\sum_{a\in\sig\cap M}\lambda_a\phi^at^a,\quad
\phi^a=\phi_1^{a_1}\dots\phi_n^{a_n}.
\end{equation}
Note that the coefficient of each monomial $t^b$ of the series \re{automorphism}
is defined by finitely many series $\lambda_a\phi^at^a$ (this is
true since the cone $\sig$ has an apex). Therefore this
is a well-defined power series. Since all the monomials in the series
belong to the semigroup $\sig\cap M$, the series defines an element of $A$.
It is easy to check that $\alpha$ is in fact a homomorphism. Also
it is clearly invertible.
\end{pf}


\subsection{Covering and combinatorial coefficients}\label{S:covering}

\begin{Def}\label{D:reasonable} Let $(X,D)$ be a toroidal pair.
We say that $(D_1,\dots,D_n)$ is a {\it reasonable covering of $D$}
if $D=D_1\cup\dots\cup D_n$, where each $D_i$ is the union of some
irreducible components of $D$, and $D_1\cap\dots\cap
D_n\subseteq\St_0(X)$.
\end{Def}

Let $(D_1,\dots,D_n)$ be a reasonable covering of $D$.
Consider an equipped local model $(X_\sig,x_0,\omega_0)$ at
a point $x\in\St_0(X)$. It can be
easily seen that the covering $(D_1,\dots,D_n)$
defines a covering of the boundary of the cone
$\sig$ in the sense of \re{covercone}, \rs{cc}.
This allows us to define the combinatorial coefficient of the
covering $(D_1,\dots,D_n)$ at each point $x\in\St_0(X)$.

\begin{Def} The {\it combinatorial
coefficient at $x\in\St_0(X)$ of the covering
$(D_1,\dots,D_n)$} is the combinatorial
coefficient of the induced covering of $\sig$ in an equipped local model
at $x$.  We denote it by $c(x)$.
\end{Def}

\begin{Rem}\label{R:invcoeff}{\bf Invariance.}
By \rr{comborder} the combinatorial coefficient
is the same for any two equipped local models that correspond
to an automorphism of $\T$ that preserves the form $\omega_0$,
and changes sign otherwise. Also it is skew-symmetric in $D_1,\dots,D_n$.
\end{Rem}

Now consider the stratification defined by the irreducible components
of $D$ (see \rs{toroidalpair}) and let
$F$ be a complete flag of stratum closures on $X$:
$$F:\ X_0\subset X_1\subset\dots\subset X_{n-1}\subset X,\quad
X_i\in\ov\St_i(X).$$
It corresponds to a complete flag of orbit closures in an equipped local
model $(X_\sig,x_0,\omega_0)$ of $(X,D)$ at $X_0$, hence to a complete
flag $F_\sig$ of faces of $\sig$.

\begin{Def} We say that the flag $F$ is positive (resp. negative) and write
$\sgn F=1$ (resp. $\sgn F=-1$) if the induced flag $F_\sig$ of faces
of $\sig$ is positive (resp. negative).
\end{Def}

Like in the case of the combinatorial coefficient the sign of
the flag depends on the choice of the form $\omega_0$ in an equipped
local model.

\begin{Def} Let $Z$ be a stratum.
We say that the closure $\overline Z$ has {\it signature}
$\{i_1,\dots,i_k\}$ for $1\leq i_l\leq n$ if and only if 
$Z\subseteq D_{i_1}\cap\dots\cap D_{i_k}$ and $k$ is maximal.
\end{Def}

The following proposition is the analog of the
description of the combinatorial coefficient given in \rt{cccone}.

\begin{Prop}\label{P:cc2}
Let $(X,D)$ be toroidal and $(D_1,\dots,D_n)$ a reasonable
covering of~$D$. Then the combinatorial coefficient $c(x)$ of
$x\in\St_0(X)$ is equal to the number of all complete flags
$$x=X_0\subset X_1\subset\dots\subset X_{n-1}\subset X,$$
where $X_i\in\ov\St_i(X)$ is a stratum closure of signature
$\{i+1,\dots,n\}$, $0\leq i\leq n-1$, counting signs.
\end{Prop}

\subsection{Symbol of monomials}

\begin{Def}\label{D:monsymbol}
Consider an ordered  collection of $n+1$ monomials in $n$ variables
with coefficients in a field $k$:
$$c_it^{a_i}=c_it_1^{a_{i1}}\dots
t_n^{a_{in}}, \quad c_i\in k^{\times}, \ \
a_i=(a_{i1},\dots,a_{in})\in\Z^n, \ \  0\leq i\leq n.$$ Let
$A=(a_{ij})\in M_{n+1,n}(\Z)$ be the matrix whose rows are the
vectors of exponents~$a_i$. Then the {\it symbol of $n+1$
monomials} is the non-zero element of $k$ defined by
$$[c_0t^{a_0},\dots,c_{n}t^{a_{n}}]=
(-1)^{B}\prod_{i=0}^{n}c_i^{(-1)^{i}A_i},$$ where $A_i$ is the
determinant of the matrix obtained from $A$ by eliminating its
$i$-th row, and
$$B=\sum_{k}\sum_{i<j}a_{ik}a_{jk}A_{ij}^k,$$ where $A_{ij}^k$ is
the determinant of the matrix obtained from $A$ by eliminating its
$i$-th and $j$-th rows and its $k$-th column.
\end{Def}

\begin{Prop}\label{P:properties} Let $f_i=c_it^{a_i}$, $0\leq i\leq n$,
be monomials. The symbol has the following properties:
\begin{enumerate}
\item (Multiplicativity) Suppose $f_i$ is a product of two
monomials $f_i=f_i'f_i''$. Then
$$[f_0,\dots,f_i'\hspace{2pt}f_i'',\dots,f_{n}]=
[f_0,\dots,f_i',\dots,f_{n}] [f_0,\dots,f_i'',\dots,f_{n}].$$
\item (Multiplicative skew-symmetry)
$$[f_0,\dots,f_{i},\dots,f_{j},\dots,f_{n}]=
[f_0,\dots,f_{j},\dots,f_{i},\dots,f_{n}]^{-1}.$$
\item (Invariance)\\
(i) Let $u=t^{Q}$, $Q\in GL(n,\Z)$ be a monomial
change of coordinates. Then
$$[\bar f_0,\dots,\bar f_{n}]=\tsymb^{\det Q},$$
where $\bar f_i=c_iu^{a_i}=c_it^{a_iQ}$ and $f_i=c_it^{a_i}$,
$0\leq i\leq n$.\\
(ii) Let $s=\lambda t$ be a translation, i.e. $s_j=\lambda_j t_j$,
$\lambda_j\in k^\times$, $1\leq j\leq n$. Then
$$[f_0',\dots,f_{n}']=\tsymb,$$
where $f_i'=c_is^{a_i}=c_i\lambda^{a_i}t^{a_i}$ and
$f_i=c_it^{a_i}$, $0\leq i\leq n$.
\end{enumerate}
\end{Prop}

\begin{pf} Modulo the sign $(-1)^B$ all the properties follow
easily from the properties of the determinant.

To take care of the sign we give an invariant description of $B$,
following \cite{Kh}. Consider $B$ as a $\Z/2\Z$-valued function of
the rows $a_0,\dots,a_{n}$ of the matrix $A$. It is easy to see
that $B=B(a_0,\dots,a_{n})$ is multilinear and its value on each
collection of $n+1$ standard vectors $(e_{i_0},\dots,e_{i_{n}})$
is 0 if more than two of the vectors $e_{i_0},\dots,e_{i_{n}}$
coincide; and 1 otherwise.

Now define the function $B'=B'(a_0,\dots,a_{n})$ to be 0 if the
rank of $(a_0,\dots,a_{n})$ is less than~$n$; and
$\lambda_0+\dots+\lambda_{n}+1\, (\,\text{mod}\, 2\,)$ if the
vectors $a_0,\dots,a_{n}$ satisfy a (unique) non-trivial
relation $\lambda_0a_0+\dots+\lambda_{n}a_{n}=0$. The function
$B'$ is multilinear and on each collection
$(e_{i_0},\dots,e_{i_{n}})$ the functions $B'$ and $B$ take the
same value. Therefore $B=B'$, in particular, $B$ is symmetric and
invariant under non-degenerate transformations.
\end{pf}


\subsection{Toroidal symbol}

Let $(X,D)$ be toroidal. Let $k(X,D)$ denote the set of rational
functions on $X$ whose divisor lies in~$D$.

Consider a 0-dimensional stratum $x\in\St_0(X)$. Then the image
of $f\in k(X,D)$ in an equipped local model
$(X_\sig,x_0,\omega_0)$ at $x$ is the product of a monomial $ct^a$
and a regular invertible function
$\phi\in\widehat\cO^\times_{X_\sig,x_0}$ with $\phi(x_0)=1$. We
call this monomial the {\it leading monomial of $f$ at $x$}. The
leading monomial is defined up to monomial transformations.

\begin{Def} Let $(X,D)$ be toroidal and $x\in\St_0(X)$ a 0-dimensional
stratum.
Define the {\it toroidal symbol $\tsymb_x$ at $x$} of
$f_0,\dots,f_{n}\in k(X,D)$ to be the symbol of the leading
monomials of $f_0,\dots,f_{n}$ at $x$.
\end{Def}

\begin{Rem}\label{R:invsymbol}{\bf Invariance.}
Let $(X_{\sig'},x_0',\omega_0')$ and
$(X_{\sig''},x_0'',\omega_0'')$ be two equipped local models at
$x$. Let $f'$ and $f''$ be the images of $f\in k(X,D)$ in the two
equipped local models. Then, according to \rp{localmodels}, the
leading monomials of $f'$ and $f''$ are related by a composition
of a monomial transformation and a translation: $t\mapsto\lambda
t^{Q}$. Therefore, by \rp{properties} the toroidal symbol is the same
for the two equipped local models if $\det Q=1$, and is reciprocal
otherwise.
\end{Rem}

By \rp{properties} the toroidal symbol is multiplicative and
multiplicatively skew-symmetric in $f_0,\dots,f_{n}$.

\medskip

Now we will give a relation between the toroidal symbol and Parshin's tame
symbol at a complete flag of stratum closures on $X$.

\begin{Prop}\label{P:relation1}
Let $(X,D)$ be toroidal. Consider $n+1$ rational functions
$f_0,\dots,f_{n}\in k(X,D)$. Then, for any complete flag $F:\
X_0\subset X_1\subset\dots\subset X_{n-1}\subset X$ of stratum
closures on $X$ we have
$$\symb_F=\tsymb_{X_0}^{\sgn F},$$
where $\symb_F$ denotes Parshin's tame symbol at the flag $F$.
\end{Prop}

(Note that the number $\tsymb_{X_0}^{\sgn F}$
is already independent of the choice of $\omega_0$ in an equipped
local model.)

\begin{pf} Since the definitions of the toroidal symbol and Parshin's
tame symbol are local we can pass to an equipped local model
$(X_\sig,x_0,\omega_0)$ at $X_0$ and assume that $F$ is a complete
flag of orbit closures in $X_\sig$.

Let $F_\sig$ be the complete flag of faces of $\sig$ corresponding to $F$,
$$F_\sig:\quad 0=\sig_0\subset\sig_1\subset\dots
\subset\sig_{n-1}\subset\sig_n=\sig.$$ Fix coordinates
$(t_1,\dots,t_n)$ in $\T$, $M\cong\Z^n$. Inside each $\sig_i$,
$1\leq i\leq n$ choose a lattice point $q_i\in\Z^n$ at
lattice distance one from $\sig_{i-1}$. Let $u_i=t^{q_i}$ be a
monomial change of coordinates in $\T$. Then the rational
functions $u_i=t^{q_i}$ give a system of local
parameters at $F$ (see appendix). Therefore by
\rp{properties} and \rr{parshsymb} of the appendix
$$\symb_F=[c_0u^{k_0},\dots,c_{n}u^{k_{n}}]=
[c_0t^{k_0},\dots,c_{n}t^{k_{n}}]^{\det Q},$$ where
$Q=(q_1,\dots,q_n)$ and $c_it^{k_i}$ is the leading monomial of
$f_i$. It remains to note that  $\det Q=\sgn F$.
\end{pf}

\begin{Cor}\label{C:relation2}
Let $(X,D)$ be toroidal and $(D_1,\dots,D_n)$ a reasonable
covering of~$D$. For $x\in\St_0(X)$ let $\cF(x)$ be the set of all
complete flags
$$x=X_0\subset X_1\subset\dots\subset X_{n-1}\subset X,$$
where $X_i\in\ov\St_i(X)$ is a stratum closure of signature
$\{i+1,\dots,n\}$, $0\leq i\leq n-1$.

Then for any $n+1$ rational functions
$f_0,\dots,f_{n}\in k(X,D)$ we have
$$\tsymb_x^{c(x)}=\prod_{F\in\cF(x)}\symb_F,$$
where we assume that the product is~1 if $\cF(x)$ is empty.
\end{Cor}

\begin{pf} This follows from \rp{cc2} and \rp{relation1}.
\end{pf}

\begin{Rem} Note that for a toroidal pair $(X,D)$ all Parshin's
reciprocity laws (\rt{Parshinrecip}) for $i>0$ follow from
\rp{relation1}. Indeed, consider a complete flag $F:\
X_0\subset\dots\subset X_i\subset\dots\subset X_n$ of irreducible
subvarieties of $X$. We can pass to a local model at $X_0$ and
assume that $X$ is an affine toric variety and $D=X\setminus\T$.
Then for any $n+1$ rational functions $f_0,\dots,f_{n}\in
k(X,D)$ the symbol $\symb_F$ is trivial unless $F$ is a flag of
orbit closures on $X$. But if we fix all $X_j$, $j\neq i$, and
vary $X_i$ there are only two such flags $F$ (since for any face
$\tau$ of a polyhedral cone there are only two codimension 1 faces
of $\tau$ that contain a fixed codimension 2 face of $\tau$), and
the signs of these two flags are opposite. Now we can apply
\rp{relation1}.
\end{Rem}


\subsection{Main theorem}

\begin{Th}\label{T:mainsymb}
Let $X$ be a complete normal $n$-dimensional variety over an
algebraically closed field $k$, and $D$ a closed subset of $X$
such that the pair $(X,D)$ is toroidal.

Let $(D_1,\dots,D_n)$ be a reasonable covering of $D$ such that each
$D_i$ is a disjoint union of two closed subsets of pure
codimension 1:
\begin{equation}\label{e:covering}
D=D_1\cup\dots\cup D_n,\ \ D_1\cap\dots\cap D_n\subseteq\St_0(X),\ \ \
D_i=D_i'\sqcup D_i'',\ \ 1\leq i\leq n.
\end{equation}
We get $2^n$ disjoint finite closed subsets of $X$:
$$S_k=G_1\cap\dots\cap G_n,\ \text{ where }
\ G_i=D_i'\text{ or }D_i'',\ \ 1\leq i\leq n,\ \ \ 1\leq k\leq
2^n.$$

Then for any $n+1$ rational
functions $f_0,\dots,f_{n}\in k(X,D)$ the following $2^n$ numbers are equal:
\begin{equation}\label{e:mainsymb}
{\Big(\prod_{x\in S_1}\tsymb_x^{c(x)}\Big)}^{(-1)^{|S_1|}}=
\dots={\Big(\prod_{x\in
S_{2^n}}\tsymb_x^{c(x)}\Big)_,}^{(-1)^{|S_{2^n}|}}
\end{equation}
where $\tsymb_x$ is the toroidal symbol of
$f_0,\dots,f_{n}$ at $x$, $c(x)$ is the combinatorial
coefficient at $x$, and $|S_k|$ is the number of $D_i''$ in the
definition of $S_k$.
\end{Th}

\begin{pf}
Because of the symmetry it is sufficient to prove the equality for
any two sets
$$S_1=D_1'\cap\dots\cap D_i'\cap\dots\cap D_n'\
\text{ and } \ S_2=D_1'\cap\dots\cap D_i''\cap\dots\cap D_n'.$$
Since the number $\tsymb_x^{c(x)}$ is multiplicatively skew-symmetric
in $D_1,\dots,D_n$ (see \rr{invcoeff}) we may assume that~$i=1$, so
$$S_1=D_1'\cap D_2'\dots\cap D_n'\
\text{ and }
\ S_2=D_1''\cap D_2'\dots\cap D_n'.$$
We have to show that
$$\prod_{x\in S_1\cup S_2}\tsymb_x^{c(x)}=1.$$

Let $\Sig$ be the union of all stratum closures $Y\in\ov\St_1(X)$
with signature $\{2',\dots,n'\}$. It follows from \rco{relation2}
that if $x\in S_1\cup S_2$ does not lie on any component of $\Sig$
then $\tsymb_x^{c(x)}=1$. On the other hand,
by \re{covering} the signature of every point
$x\in\St_0(X)\cap\Sig$ is either $$\{2',\dots,n'\},\text{ or }
\{1',2',\dots,n'\},\text{ or }\{1'',2',\dots,n'\}.$$ In the first
case $\tsymb_x^{c(x)}=1$, again by \rco{relation2}. In the second
case $x\in S_1$ and in the third $x\in S_2$. Therefore, we have
\begin{equation}\label{e:1}
\prod_{x\in S_1\cup S_2}\tsymb_x^{c(x)}\ =
\prod_{x\in\St_0(X)\cap\Sig}\tsymb_x^{c(x)}.
\end{equation}

Now consider a component $Y$ of $\Sig$, and a closed point~$y\in Y$. Let
$\cF(y,Y)$ denote the set of all complete flags
$$y\subset Y\subset X_2\subset\dots\subset X_{n-1}\subset X,$$
where $X_i\in\ov\St_i(X)$ has signature  $\{(i+1)',\dots,n'\}$,
$2\leq i\leq n-1$.
Denote $$\symb_{y,Y}=\prod_{F\in\cF(y,Y)}\symb_{F},$$
and we assume that $\symb_{y,Y}=1$ if $\cF(y,Y)$ is empty. Then by
\rco{relation2} for each $x\in\St_0(X)\cap\Sig$ we have
\begin{equation}\label{e:2}
\tsymb_x^{c(x)}=\prod_{F\in\cF(x)}\symb_F=\prod_{Y\ni
x}\symb_{x,Y},
\end{equation}
where the product on the right hand side runs over all components
$Y$ of $\Sig$ containing~$x$.

On the other hand, by the first Parshin's reciprocity law
(\rt{Parshinrecip} for $i=0$)
$$\prod_{y\in Y}\symb_{y\subset Y\subset X_2\subset\dots\subset X_{n-1}\subset X}=1,$$
where the product is taken over all points $y\in Y$. Thus
$$\prod_{y\in Y}\symb_{y,Y}=1.$$
Note that $\symb_{y,Y}$ is trivial for all points $y$ not lying
in~$\St_0(X)$, so we can assume that $y\in\St_0(X)\cap Y$. We have
\begin{equation}\label{e:3}
\prod_{y\in\St_0(X)\cap Y}\symb_{y,Y}=1.
\end{equation}

Combining \re{1}, \re{2} and \re{3} we get
\begin{eqnarray}
\prod_{x\in S_1\cup S_2}\tsymb_x^{c(x)}&=&
\prod_{x\in\St_0(X)\cap\Sig}\tsymb_x^{c(x)}=\nonumber\\
\prod_{x\in\St_0(X)\cap\Sig}\ \ \prod_{Y\ni x}\symb_{x,Y}&=&
\prod_{Y\subset\Sig}\ \ \prod_{x\in\St_0(X)\cap
Y}\symb_{x,Y}=1.\nonumber
\end{eqnarray}
\end{pf}


\section{Toroidal Residue}\label{S:residue}

Let $(X,D)$ be a toroidal pair, as before. At each point
$x\in\St_0(X)$ we define the residue $\tres_x\omega$ of a rational
$n$-form $\omega$ on $X$ which is regular in $X\setminus D$. Then we prove
an additive analog of \rt{mainsymb}.

\subsection{Toroidal residue}
First we will define the toroidal residue for a local model
$(X_\sig,x_0,\omega_0)$ at a point $x\in\St_0(X)$. As before
$X_\sig$ is an $n$-dimensional affine toric variety, $x_0$ is the closed
orbit, and
$\omega_0=\frac{dt_1}{t_1}\wedge\dots\wedge\frac{dt_n}{t_n}$,
where $(t_1,\dots,t_n)$ is a coordinate system in $\T$.

Let $A=\hat\cO_{X_\sig,x_0}$ be the completion of the local ring
of $x_0$ on $X_\sig$, $B=A_S$ the localization of $A$ by the
multiplicative subgroup $S$ of all monomials. We consider the
$B$-algebra $\Omega_{B}^n$ of differential $n$-forms that are
regular in $\T$, and the $A$-algebra $\Omega_{A}^n$ of regular
differential $n$-forms.

By fixing coordinates $(t_1,\dots,t_n)$ we can identify
every element $f\in B$ with a formal power series
$$f=t^b\sum_{a\in\sig\cap\Z^n}\lambda_at^a,\quad b\in\Z^n.$$
Let $\omega\in\Omega_{B}^n$. Then we can write
$\omega=f\frac{dt_1}{t_1}\wedge\dots\wedge\frac{dt_n}{t_n}$,
for some $f\in B$.

\begin{Def} The {\it toroidal residue} of a differential $n$-form
$\omega\in\Omega_B^n$ is the constant term $\lambda_{-b}$ in the
formal power series of $f$. We denote it by $\tres\omega$.
\end{Def}

We have the following properties of the toroidal residue.
\begin{Prop}\label{P:resproperties} Consider a differential
$n$-form $\omega\in\Omega_B^n$. Then
\begin{enumerate}
\item If $\omega$ is exact then $\tres\omega=0$;
\item If $\omega\in\Omega_{A}^n$ then $\tres\omega=0$;
\item For any $s_1,\dots,s_n\in B^{\times}$,
$\tres\frac{ds_1}{s_1^{m_1}}\wedge\dots\wedge\frac{ds_n}{s_n^{m_n}}=0$,
unless all $m_i=1$;
\item The toroidal residue is independent of the choice of coordinates
$(t_1,\dots,t_n)$;
\item The toroidal residue is invariant under monomial transformations
$t\mapsto t^Q$, $Q\in GL(n,\Z)$ up to a factor $\det Q$.
\end{enumerate}
\end{Prop}

\begin{pf}
i) Let $\omega=d{w}$. We can assume that
${w}=g{dt_1}{}\wedge\dots\wedge{\widehat{dt_i}}{}
\wedge\dots\wedge{dt_n}{}$. Then
$$\omega=(-1)^{i-1}\frac{\partial g}{\partial t_i}t_1\dots t_n
\frac{dt_1}{t_1}\wedge\dots\wedge\frac{dt_n}{t_n}.$$
Suppose $g=\sum_a\lambda_at^a$. Then
$$(-1)^{i-1}\frac{\partial g}{\partial t_i}t_1\dots t_n=
\sum_{a}a_i\lambda_a\frac{t^a}{t_i}t_1\dots t_n,\ \ a=(a_1,\dots,a_n).$$
Clearly the constant term of the last series is zero.\newline

\hspace{-0.42cm}ii) Let $u_i=t^{a_i}$, $a_i\in\sig\cap\Z^n$ be $n$
regular functions, such that $du_1,\dots,du_n$ are linearly
independent. Then for every $\omega\in\Omega_{A}^n$
$$\omega=fdu_1\wedge\dots\wedge du_n=fdt^{a_1}\wedge\dots\wedge dt^{a_n}=
ft^{a_1+\dots+a_n}J\frac{dt_1}{t_1}\wedge\dots\wedge\frac{dt_n}{t_n},$$
where $f\in A$ and $J=\det(a_1,\dots,a_n)$.
Clearly, the constant term of $ft^{a_1+\dots+a_n}J$
is zero.\newline

\hspace{-0.42cm}iii) First assume that $char(k)=0$. Suppose
$m_i\neq 1$. Then
$$\frac{ds_1}{s_1^{m_1}}\wedge\dots\wedge\frac{ds_n}{s_n^{m_n}}=
d\left(\frac{(-1)^{i-1}}{1-m_i}s_i^{1-m_i}\frac{ds_1}{s_1^{m_1}}\wedge\dots
\wedge\frac{\widehat{ds_i}}{s_i^{m_i}}\wedge\dots\wedge
\frac{ds_n}{s_n^{m_n}}\right),$$ and the statement follows from
part~i). In the case of an arbitrary characteristic note that
the toroidal residue is a polynomial function in finitely many
coefficients of the series $s_1,\dots,s_n\in B^{\times}$. This
function is independent of the characteristic and vanishes when
the characteristic is zero. Therefore it is identically
zero.\newline

\hspace{-0.42cm}iv) Let $(s_1,\dots,s_n)$ be another coordinate
system in $\T$. Then $s_i=\phi_it_i$, where $\phi_i\in A^\times$.
Consider $\omega\in\Omega_B^n$ and write
$$\omega=f\frac{ds_1}{s_1}\wedge\dots\wedge\frac{ds_n}{s_n},\quad
f=\sum_a\lambda_as^a,\ a\in\Z^n.$$
Then the residue of $\omega$ with respect to $(s_1,\dots,s_n)$ is
$\tres_{(s_1,\dots,s_n)}\omega=\lambda_0$. On the other hand, by part~iii)
the residue of $\omega$ with respect to $(t_1,\dots,t_n)$ equals
\begin{equation}\label{e:open}
\tres_{(t_1,\dots,t_n)}\omega=
\tres_{(t_1,\dots,t_n)}
\lambda_0\frac{ds_1}{s_1}\wedge\dots\wedge\frac{ds_n}{s_n}.
\end{equation}
Now for each $1\leq i\leq n$,
$\frac{ds_i}{s_i}=\frac{d\phi_i}{\phi_i}+\frac{dt_i}{t_i}$. Substituting
into \re{open} and expanding we get
$$\tres_{(t_1,\dots,t_n)}\omega=\tres_{(t_1,\dots,t_n)}\lambda_0
\frac{dt_1}{t_1}\wedge\dots\wedge\frac{dt_n}{t_n}+
\sum_i\tres_{(t_1,\dots,t_n)}\omega_i,$$ where in each $\omega_i$
at least one of $\frac{dt_j}{t_j}$ is replaced by
$\frac{d\phi_j}{\phi_j}$. It is easy to check that the residue of
every $\omega_i$ is zero.
\newline

\hspace{-0.42cm}v) It follows from the fact that if
$u_i=t^{q_i}$, $q_i\in\Z^n$ then
$$\frac{du_1}{u_1}\wedge\dots\wedge\frac{du_n}{u_n}=\det Q
\frac{dt_1}{t_1}\wedge\dots\wedge\frac{dt_n}{t_n},$$ and the
observation that monomial transformations do not change the
constant term of a series.
\end{pf}

\begin{Rem} The proof of parts iii) and iv) is similar
to the one given in \cite{F-P} for Parshin's residue.
\end{Rem}

Now let $(X,D)$ be toroidal. Denote by $\Omega^n(X,D)$ the
set of all rational $n$-forms on $X$ that are regular in $X\setminus D$.

\begin{Def} The {\it toroidal residue} of a rational $n$-form
$\omega\in\Omega^n(X,D)$ at a point $x\in\St_0(X)$ is the toroidal residue
of its image in a local model at $x$. We denote it by $\tres_x\omega$.
\end{Def}

\begin{Rem}\label{R:invres}{\bf Invariance.}
As it follows from \rp{localmodels} and \rp{resproperties}
the toroidal residue is the same for any two equipped local models that
are related by an isomorphism preserving the form $\omega_0$, and
changes sign otherwise.
\end{Rem}

The relation between the toroidal residue and Parshin's residue is
similar to the one between the toroidal symbol and Parshin's tame symbol.

\begin{Prop}\label{P:resrelation1}
Let $(X,D)$ be toroidal. Consider a rational $n$-form $\omega\in\Omega^n(X,D)$.
Then for any complete flag
$F:\ X_0\subset X_1\subset\dots\subset X_{n-1}\subset X$ of
stratum closures on $X$ we have
$$\res_F\omega=\sgn F\tres_{X_0}\omega,$$
where $\res_F\omega$ denotes Parshin's residue at the flag $F$.
\end{Prop}

The number $\sgn F\tres\omega$ is independent of the choice of the
form $\omega_0$ in an equipped local model. Consequently, if
$(D_1,\dots,D_n)$ is a reasonable covering of $D$ then at each
$x\in\St_0(X)$ the number $c(x)\tres_x\omega$ is also
well-defined.

The following is the additive analog of \rco{relation2}.

\begin{Cor}\label{C:resrelation2}
Let $(X,D)$ be toroidal and $(D_1,\dots,D_n)$ a reasonable
covering of~$D$. For $x\in\St_0(X)$ let $\cF(x)$ be the set of all
complete flags
$$x=X_0\subset X_1\subset\dots\subset X_{n-1}\subset X,$$
where $X_i\subset\ov\St_i(X)$ is a stratum closure of signature
$\{i+1,\dots,n\}$, $0\leq i\leq n-1$.

Then for any $\omega\in\Omega^n(X,D)$
we have
$$c(x)\tres_x\omega=\sum_{F\in\cF(x)}\res_F\omega,$$
where we assume that the sum is~0 if $\cF(x)$ is empty.
\end{Cor}


\subsection{Main theorem}
\begin{Th}\label{T:mainres} Let $X$ be a complete normal
$n$-dimensional variety over an algebraically closed field $k$,
and $D$ a closed subset of $X$ such that the pair $(X,D)$ is
toroidal.

Let $(D_1,\dots,D_n)$ be a reasonable covering of $D$ such that
each $D_i$ is a disjoint union of two closed subsets of pure
codimension 1:
\begin{equation}\label{e:rescover}
D=D_1\cup\dots\cup D_n,\ \ D_1\cap\dots\cap D_n\subset\St_0(X),\ \
\ D_i=D_i'\sqcup D_i'',\ \ 1\leq i\leq n.
\end{equation}
We get $2^n$ disjoint finite closed subsets of $X$:
$$S_k=G_1\cap\dots\cap G_n,\ \text{ where }
\ G_i=D_i'\text{ or }D_i'',\ \ 1\leq i\leq n,\ \ \ 1\leq k\leq
2^n.$$

Then for any rational $n$-form $\omega\in\Omega^n(X,D)$,
the following $2^n$ numbers are equal:
$$(-1)^{|S_1|}\sum_{x\in S_1}c(x)\tres_x\omega=\dots=
(-1)^{|S_{2^n}|}\sum_{x\in S_{2^n}}c(x)\tres_x\omega,$$ where
$\tres_x\omega$ denotes the toroidal residue of $\omega$ at $x$,
$c(x)$ is the combinatorial coefficient at $x$, and $|S_k|$ is the
number of $D_i''$ in the definition of $S_k$.
\end{Th}

\begin{pf} The proof repeats the arguments of the proof
of \rt{mainsymb}.
\end{pf}


\section{Application to Systems of Equations in the Torus}\label{S:application}

In this section we apply our main results on the toroidal symbol and
residue to prove the product of roots formula and the sum of values formula
(Theorems A and B in the introduction).

Recall that a {\it Laurent polynomial} is a finite linear combination of
monomials with integer exponent vectors and coefficients in $k$:
$$f(t)=\sum_{m\in\Z^n}\lambda_mt^m,\quad t^m=t_1^{m_1}\dots t_n^{m_n},\ \ \lambda_m\in k.$$
The convex hull of those lattice points $m\in\Z^n$ for which $\lambda_m\neq 0$ is
called the {\it Newton polytope of $f$}. The value of $f(t)$ is defined for all $t$
in the algebraic $n$-torus $\T=(k^\times)^n$.

Consider a system of $n$ Laurent polynomial equations in $\T$:
\begin{equation}\label{e:intropoly}
f_1(t)=\dots=f_n(t)=0,\quad t\in\T.
\end{equation}
Let $\D_i$ be the Newton polytope of $f_i$. 
We assume that none of the $f_i$ is a monomial, hence,
none of $\D_i$ is a point.

Every linear functional ${w}$ on $\R^n$ defines a collection of faces
$\D_1^{w},\dots,\D_n^{w}$ of the Newton polytopes such that the
restriction of ${w}$ on $\D_i$ attains its maximum precisely at
$\D_i^{w}$. The polynomial 
$$f_i^{w}(t)=\sum_{m\in\D_i^{{w}}\cap\Z^n}\lambda_mt^m$$
is called the {\it initial form} of $f_i$ with respect to ${w}$.
According to Bernstein's theorem \cite{B} the number of the solutions
to the system \re{intropoly} is finite (and equals the mixed volume of $\D_1,\dots,\D_n$)
if and only if for every ${w}\neq 0$ the system 
$f_1^{w}(t)=\dots=f_n^{w}(t)=0$ is
inconsistent.

\begin{Def}\label{D:developed}
A collection of polytopes $\D_1,\dots,\D_n$ is called {\it developed} if
none of them is a point and for each ${w}\neq 0$ at least one
of the faces $\D_1^{w},\dots,\D_n^{w}$ is a vertex.
\end{Def}

By above any system with developed collection of Newton polytopes
has finitely many solutions. We will call them the {\it roots} of the system.
Every root $x$ of the system has a multiplicity $\mu(x)$.

Let $\D$ be the Minkowski sum of $\D_1,\dots,\D_n$. Then
every face $\G\subset\D$ has a unique decomposition as a
sum of faces
\begin{equation}\label{e:introdecomp}
\G=\G_1+\dots+\G_n,\ \ \ \text{ where } \G_i\subset\D_i,\ \
i=1,\dots,n.
\end{equation}

If the collection $\D_1,\dots,\D_n$ is developed then $\D$ has
dimension $n$ and in the decomposition of every proper face of
$\D$ at least one summand is a vertex. In this case, for each
vertex $A$ of $\D$ the combinatorial coefficient $c(A)$ is defined.

\begin{Def} Let $\sig_A$ be the cone with apex $A$ generated by the
facets of $\D$ that contain~$A$. Then the boundary of $\sig_A$ is
covered by the closed sets $D_1,\dots,D_n$, where $D_i$ is the
union of all facets of $\sig_A$ that correspond to the facets of $\D$
whose $i$-th summand in the decomposition \re{introdecomp} is a
vertex. The combinatorial coefficient of this covering is called
the {\it combinatorial coefficient} $c(A)$ of the vertex $A\in\D$.
\end{Def}


\subsection{Product of roots formula}
Consider a system of $n$ Laurent polynomials
with developed collection of Newton polytopes.
The product of the roots counting multiplicities is a
point in $\T$, which we denote by $\rho$. To locate $\rho$ it is enough to
find the product of the values of $t_i$ 
over the roots of the system, for each $1\leq i\leq n$.
More generally, we will find the product of the values of 
any Laurent monomial $ct^m$, for $c\in k^\times$, $m\in\Z^n$, 
over the roots of the system. 

\begin{Def}\cite{Kh}
The {\it symbol of $f_1,\dots,f_n$ and a Laurent monomial $f_0=ct^m$ 
at a vertex $A\in\D$} is the symbol of $n+1$ monomials
$[ct^m,f_1(A_1)t^{A_1},\dots,f_n(A_n)t^{A_n}]$, where
$A=A_1+\dots+A_n$ is the decomposition of $A$, and $f_i(A_i)$ is
the coefficient of $t^{A_i}$ in $f_i$. We denote it by
$[f_0,\dots,f_n]_A$.
\end{Def}

The following theorem was proved by Khovanskii \cite{Kh}
in the complex case. Our proof uses the result of \rt{mainsymb}
and works over an arbitrary algebraically closed field $k$.

\begin{Th}\label{T:product}
Suppose the collection of the Newton polytopes $\D_1,\dots,\D_n$
of the system \re{intropoly} is developed. Then the product of the
values of a Laurent monomial $f_0$ over the roots
of the system is given by
\begin{equation}\label{e:prod}
\prod_{f_1(x)=\dots=f_n(x)=0}f_0(x)^{\mu(x)}=
\prod_{A}[f_0,\dots,f_n]_A^{(-1)^nc(A)},
\end{equation}
where the right hand product
is taken over all vertices $A$ of the polytope
$\D=\D_1+\dots+\D_n$ and $c(A)$ is the combinatorial coefficient
at $A$.
\end{Th}

\begin{pf} First let us notice that it is sufficient to prove the theorem for
a generic system with given Newton polytopes. There is an open set
$U$ in the space of the coefficients of the system where the
number of the roots counting multiplicities is constant. The
left hand side of \re{prod} being symmetric in
the roots of the system is a rational function in the
coefficients of the system. On the other hand, the product of the
symbols $[f_0,\dots,f_n]_A$ is also a rational function of
the coefficients of the system. Suppose we proved the formula for
almost all systems in $U$. Then the two rational functions
coincide on an open algebraic subset $W\subset U$, thus coincide
everywhere in $U$.

Consider the complete toric variety $X$ associated with the Minkowski sum
$\D$ (e.g. see \cite{F}).
Let $D=X\setminus\T$ be the invariant divisor on $X$.
Denote by $Z_i$ the closure of the
zero locus $f_i=0$ in $X$, and let $Z=Z_1\cup\dots\cup Z_n$. If the
system is generic, the components of $Z$
intersect transversally and the intersection of each component of
$Z$ with $D$ is also transversal. In this case the pair $(X,D\cup
Z)$ is toroidal.

Now we will define a covering of $D\cup Z$. Each irreducible
component of $D$ corresponds to a facet of $\D$. Recall that each
facet $\G$ of $\D$ has a unique decomposition into the sum of
faces \re{introdecomp}. Denote by $D_i$ the union of all
components that correspond to those facets whose $i$-th summand is
a vertex. Since the collection $\D_1,\dots,\D_n$ is developed, the
sets $D_1,\dots,D_n$ define a covering of $D$.

Consider a covering $D\cup Z=(D_1\cup Z_1)\cup\dots\cup(D_n\cup
Z_n)$. Notice that every intersection point in $(D_1\cup Z_1)\cap\dots\cap(D_n\cup
Z_n)$ is either a fixed orbit or a transversal intersection of some components
of $Z$ and $D$. Thus the covering is reasonable (see \rs{covering}).
By definition, the components of $D_i$ correspond to facets of $\D$ whose $i$-th
summand is a vertex, i.e. the corresponding initial forms of $f_i$ are monomials.
This implies that $D_i\cap Z_i=\varnothing$. Now applying
\rt{mainsymb} we obtain
$$\prod_{x\in Z_1\cap\dots\cap Z_n}[f_0,\dots,f_n]_x^{c(x)}=
\prod_{x\in D_1\cap\dots\cap D_n}[f_0,\dots,f_n]_x^{(-1)^nc(x)}.$$

It remains to notice that for each transversal intersection $x$ of
components of $Z$ the combinatorial coefficient $c(x)=1$ and
$[f_0,\dots,f_n]_x=f_0(x)^{\mu(x)}$ (see \rex{example} in the
appendix). Also for a point $x\in D_1\cap\dots\cap D_n$ the toric
symbol $[f_0,\dots,f_n]_x$ coincides with
$[f_0,\dots,f_n]_A$, where $A$ is the corresponding vertex of
$\D$, and $c(x)=c(A)$, by definition.
\end{pf}


\subsection{Sum of values formula}

We recall the definition of a Laurent series at a vertex and
the residue at a vertex from \cite{G-Kh}.
Let $f$ be a Laurent polynomial with Newton polytope $\D(f)$,
$A$ a vertex of $\D(f)$, and $f(A)$ is the coefficient of 
$t^A$ in $f$. Since the constant term of the Laurent
polynomial $\tilde f=f/(f(A)t^A)$ equals~1 we get a
well-defined power series
\begin{equation}\label{e:introseries}
\frac{1}{\tilde f}=1+(1-\tilde f)+(1-\tilde f)^2+\dots.
\end{equation}

\begin{Def} Let $g$ be a Laurent polynomial. The formal product of
the series \re{introseries} and the Laurent
polynomial $g/(f(A)t^A)$ is called the {\it Laurent series of
$g/f$ at the vertex $A\in\D(f)$}.
\end{Def}

\begin{Def} The {\it residue at a vertex $A\in\D(f)$} of a
rational $n$-form
$\omega=
\frac{g}{f}\,\frac{dt_1}{t_1}\wedge\dots\wedge\frac{dt_n}{t_n}$
is the constant term of the Laurent series of $g/f$ at $A$. We denote
it by  $\res_A\omega$.
\end{Def}

The following theorem was proved by Gel'fond and Khovanskii \cite{G-Kh2}
in the case when $k=\C$.
We prove it for any algebraically closed field $k$ using \rt{mainres}.

\begin{Th}\label{C:introsum} Suppose the collection of the Newton
polytopes $\D_1,\dots,\D_n$ of the system \re{intropoly} is
developed. Then the sum of the values of a Laurent polynomial $f_0$
over the roots of the system counting multiplicities is given by
$$\sum_{f_1(x)=\dots=f_n(x)=0}\mu(x)f_0(x)=(-1)^n\sum_{A}c(A)\res_A
\left(\frac{f_0\,J}{f}\,
\frac{dt_1}{t_1}\wedge\dots\wedge\frac{dt_n}{t_n}\right),$$ where
the sum on the right is taken over the vertices $A$ of
$\D=\D_1+\dots+\D_n$, 
$J=\det\left({t_j}\frac{\partial f_i}{\partial t_j}\right)$ for 
$1\leq i,j\leq n$, $f=f_1\dots f_n$,
and $c(A)$ is the combinatorial coefficient at $A$.
\end{Th}

\begin{pf} As in the proof of \rt{product} it is enough to consider the case of
a generic system with given Newton polytopes, since the sum of
the values of a Laurent polynomial over the roots of the system is
a rational function of the coefficients of the system.

As before, $X$ is the complete toric variety associated with the
Minkowski sum, $D$ the invariant divisor on $X$,
$Z_i$ the closure of the zero locus $f_i=0$ in $X$, and $D\cup
Z=(D_1\cup Z_1)\cup\dots\cup(D_n\cup Z_n)$ a reasonable covering
that satisfies $D_i\cap Z_i=\varnothing$. Applying \rt{mainres} to
the form $\omega=f_0\frac{df_1}{f_1}\wedge\dots\wedge\frac{df_n}{f_n}$ 
we get
$$\sum_{x\in Z_1\cap\dots\cap Z_n}c(x)\tres_x\omega=(-1)^n
\sum_{x\in D_1\cap\dots\cap D_n}c(x)\tres_x\omega.$$

For each transversal intersection $x$ of
components of $Z$ the combinatorial coefficient $c(x)=1$ and
$\tres_x\omega=\mu(x)f_0(x)$ (see \rex{example2} in the appendix). Also
for a point $x\in D_1\cap\dots\cap D_n$, we have $c(x)=c(A)$ and
$$\tres_x\omega=\res_A
\left(\frac{f_0\,J}{f}\,
\frac{dt_1}{t_1}\wedge\dots\wedge\frac{dt_n}{t_n}\right),$$ where
$A\in\D$ is the vertex corresponding to the fixed orbit $x$.
Therefore, we have obtained the required equality.
\end{pf}


\section{Appendix: Parshin's Reciprocity Laws}

Here we recall the definition of Parshin's tame symbol and residue
for an arbitrary algebraic variety $X$ over an algebraically
closed field, and formulate Parshin's reciprocity laws.

\subsection{Parshin's tame symbol}

Let $X$ be a complete algebraic variety over an algebraically
closed field $k$.

Consider a complete flag of irreducible subvarieties of $X$:
\begin{equation}\label{e:pflag}
F:\ \ X_0\subset X_{1}\subset\dots\subset X_{n-1}\subset X_n=X.
\end{equation}
We will assume that all $X_i$ are normal. The general case
can be reduced to this one by considering normalization (for
details see \cite{P2} or \cite{T}). Note also that this
assumption holds for complete flags of stratum closures on
toroidal pairs.

Given a flag $F$ as in \re{pflag} define a
{\it system of local parameters at $F$}, $(u_1,\dots,u_n)_F\in k(X)^n$ 
as follows.
There exists an open subset of $X_n$ where the
codimension 1 subvariety $X_{n-1}$ has a local equation
$u_n$. In general, for every $i=0\dots n-1$, let $u_{n-i}$ be a
local equation (in some open subset) of the codimension~1
subvariety $X_{n-i-1}\subset X_{n-i}$.

Next, for every rational function $f$ on $X$ define its {\it
order} $(a_1,\dots,a_n)_F\in\Z^n$ at the flag $F$.  
First let $a_n$ be the order
of $f$ along $X_{n-1}$. We
can write
$$f=f^{(n-1)}u_n^{a_n},\quad a_n\in\Z.$$
Let $\tilde f^{(n-1)}$ be the restriction of $f^{(n-1)}$ to $X_{n-1}$,
and $a_{n-1}$ be the order of this restriction along
$X_{n-2}$,
$$\tilde f^{(n-1)}=f^{(n-2)}u_{n-1}^{a_{n-1}},\quad a_{n-1}\in\Z,$$
and so on. Finally,
\begin{equation}\label{e:above}
\tilde f^{(1)}=f^{(0)}u_1^{a_1},\quad a_1\in\Z,
\end{equation}
where $a_1$ is the order of $\tilde f^{(1)}$ at $X_0$.

\begin{Def} Let $f_0,\dots,f_{n}$ be rational functions
on $X$. Fix a complete flag $F$ of irreducible subvarieties
\re{pflag}. Let
$(a_{i1},\dots,a_{in})_F$ be the order of $f_i$ at
$F$. Denote $A=(a_{ij})\in M_{n+1,n}(\Z)$.
{\it Parshin's tame symbol} of $f_0,\dots,f_{n}$ at the flag $F$
is the following non-zero element of $k$:
\begin{equation}\label{e:psymbol}
\symb_{F}=(-1)^B\Big(\prod_{i=0}^{n}f_i^{(-1)^{i}A_i}\Big)(X_0),
\end{equation}
where $A_i$ is the determinant of the matrix obtained from $A$ by
eliminating its $i$-th row, and
$$B=\sum_{k}\sum_{i<j}a_{ik}a_{jk}A_{ij}^k,$$ where $A_{ij}^k$ is
the determinant of the matrix obtained from $A$ by eliminating its
$i$-th and $j$-th rows and its $k$-th column.
\end{Def}

Note that the order of the rational function inside the large
brackets in \re{psymbol} is~$(0,\dots,0)_F$, hence, its value at
$X_0$ makes sense and is not zero.

\begin{Rem}\label{R:parshsymb}
Let us associate with every rational function $f$ on $X$ a
monomial $cu_1^{a_1}\dots u_n^{a_n}$, where $(u_1,\dots,u_n)_F$ are
local parameters at $F$, $(a_1,\dots,a_n)_F$ is the order of $f$
at $F$, and $c=f^{(0)}(X_0)$. Then the tame symbol
$\symb_F$ is equal to the symbol of the
corresponding $n+1$ monomials (see \rd{monsymbol}).
\end{Rem}

Parshin's symbol does not depend on the
choice of local parameters $(u_1,\dots,u_n)_F$. It is
multiplicative and skew-symmetric (compare to~\rp{properties}).

\begin{Example}\label{ex:example}
Let $f_1,\dots,f_n$ be rational functions on an algebraic variety
$X$, whose zero loci $\{f_i=0\}$ intersect transversely at a
non-singular point $x\in X$. Denote by $Z_i$ the irreducible
component of $\{f_i=0\}$ that contains $x$.

Let $f_0$ be any rational function on $X$ whose divisor does not
contain $x$. Then
$$\langle f_0,\dots,f_n\rangle_F=f_0(x)^{\mu(x)},$$
where $F:\,x=X_0\subset X_1\subset\dots\subset X_{n-1}\subset X$ for
$X_i=Z_{i+1}\cap\dots\cap Z_{n}$, and $\mu(x)=\mu_1\dots\mu_n$ is the
product of the multiplicities $\mu_i$ of $f_i$ along $Z_i$.
\end{Example}

Indeed, for the system of local parameters at $F$ we can choose
the local equations of $Z_i$ at $x$. Then the first row
of the matrix $A$ is zero and the other $n$ rows form a lower
triangular matrix with the multiplicities $\mu_i$ on the diagonal.
Therefore, $A_{0}=\mu_1\dots\mu_n$ and $A_j=0$, $1\leq j\leq n$.
It is also not hard to see that $B=0$ (for example one can use the
description of $B$ given in the proof of \rp{properties}).


\subsection{Parshin's residue}

Let $X$ be a complete algebraic variety over an
algebraically closed field $k$, and $F$ a complete flag of
irreducible subvarieties \re{pflag}.
Let $(u_1,\dots,u_n)_F$ be a system
of local parameters at $F$, as before.

Consider a rational differential $n$-form on $X$.
At a generic point of $X_{n-1}$ the
differentials $du_1,\dots,du_n$ are linearly independent, and we
can write
$$\omega=fdu_1\wedge\dots\wedge du_n,\ \ \text{ where }\ \
f=\sum_{i_n>N_n}f_{i_n}u_n^{i_n}.$$ The restriction of the form
$f_{-1}du_1\wedge\dots\wedge u_{n-1}$ onto $X_{n-1}$ makes
sense and gives us a rational $(n-1)$-form $\omega_{n-1}$ on
$X_{n-1}$. Continuing in this way we arrive at a sequence
of rational $(n-i)$-forms $\omega_{n-i}$ on $X_{n-i}$,
$i=1,\dots,n$, the last one being a number
$\omega_0=f_{-1,\dots,-1}$ at the point $X_0$. Note also that
this number is the coefficient of the series
$$f=\sum_{i_n\geq N_n}\ \sum_{i_{n-1}\geq N_{n-1}(i_{n})}\dots
\sum_{i_1\geq N_1(i_2,\dots,i_n)}f_{i_1,\dots,i_n}u_1^{i_1}\dots
u_n^{i_n},\ \ \ f_{i_1,\dots,i_n}\in k,$$  where we identify $f$
with an element of the field $k((u_1))\dots((u_n))$. Here
$K((t))$ denotes the field of the Laurent power series in $t$ with
coefficients in a field $K$.

\begin{Def}\label{D:Parshres}
Let $\omega$ be a rational $n$-form on $X$. Fix a complete flag
$F$ of irreducible subvarieties \re{pflag}.
{\it Parshin's residue $\res_F\omega$} at the flag $F$ is
the number $f_{-1,\dots,-1}$
constructed above.
\end{Def}

Parshin's residue does not depend on the choice of local
parameters $(u_1,\dots,u_n)_F$. The proof of
this statement is similar to the proof we gave for the invariance
of the toroidal residue in \rp{resproperties}.

\begin{Example}\label{ex:example2}
Let $f_1,\dots,f_n$ be rational functions on an algebraic variety
$X$, whose zero loci $\{f_i=0\}$ intersect transversely at a
non-singular point $x\in X$. Denote by $Z_i$ the irreducible
component of $\{f_i=0\}$ that contains $x$.

Let $f_0$ be any rational function on $X$ which is regular in an open
neighborhood of $x$. Then
$$\res_Ff_0\frac{df_1}{f_1}\wedge\dots\wedge\frac{df_n}{f_n}=\mu(x)f_0(x),$$
where $F:\, x=X_0\subset X_1\subset\dots\subset X_{n-1}\subset X$, for
$X_i=Z_{i+1}\cap\dots\cap Z_{n}$, and $\mu(x)=\mu_1\dots\mu_n$ is the
product of the multiplicities $\mu_i$ of $f_i$ along $Z_i$.
\end{Example}

Indeed, for the system of local parameters at $F$ we can choose
the local equations $u_i$ of $Z_i$ at $x$. Then
$$\res_F f_0\frac{df_1}{f_1}\wedge\dots\wedge\frac{df_n}{f_n}=
\res_F f_0\frac{du_1^{\mu_1}}{u_1^{\mu_1}}\wedge\dots\wedge
\frac{du_n^{\mu_n}}{u_n^{\mu_n}}=
\res_F\mu(x)f_0\frac{du_1}{u_1}\wedge\dots\wedge\frac{du_n}{u_n}=
\mu(x)f_0(x).$$


\subsection{Reciprocity laws}

Now we will formulate Parshin's reciprocity laws for the tame
symbol and the residue.

\begin{Th}\label{T:Parshinrecip}
Let $X$ be a complete irreducible $n$-dimensional algebraic
variety over an algebraically closed field $k$. Fix a partial flag
of irreducible subvarieties $X_0\subset\dots\subset\widehat
X_i\subset\dots\subset X_n=X$, where $X_i$ is omitted. Then
\begin{enumerate}
\item for any $n+1$ rational functions $f_0,\dots,f_{n}$ on $X$
$$\prod_{X_i}\symb_{X_0\subset\dots\subset X_i\subset\dots\subset X_n}=1,$$
\item for any rational $n$-form $\omega$ on $X$
$$\sum_{X_i}\res_{X_0\subset\dots\subset X_i\subset\dots\subset X_n}\omega=0,$$
\end{enumerate}
where the product (sum) is taken over all irreducible subvarieties
$X_i$ that complete the flag, and is finite.
\end{Th}

It follows by definition that Parshin's symbol $\symb_{F}$ equals 1
unless $F$ consists of intersections of components of the divisors of $f_0,\dots,f_{n}$.
This shows that the above product is finite. Similarly, Parshin's residue
$\res_F\omega$ is zero unless $F$ consists of intersections of components of 
the polar set of $\omega$, hence, the above sum is finite.

In the proof of the main theorems (\rt{mainsymb}, \rt{mainres}) we referred to the special case
of the reciprocity law when $i=0$. In this case the proof of the
first part of \rt{Parshinrecip} is based on the ``reduction formula'' for the symbol.
Namely, the property of the symbol analogous to the cofactor
expansion for the determinant allows one to represent the
$n$-dimensional symbol as a product of symbols of dimension $n-1$
(see \cite{F-P}). Then the statement follows from Weil's reciprocity
by induction.

For the case of the residue, notice that the residue of the
$n$-form $\omega$ at the flag $F$ is equal to the sum of the
residues of $1$-forms $\omega_1$ on $X_1$ at $X_0$. The
statement then follows from the 1-dimensional residue formula.

\bibliographystyle{amsalpha}

\end{document}